\documentclass[number,citesort,seceqn,dvips]{arxbj}

\aid{0}
\volume{17}
\issue{2}
\pubyear{2011}
\firstpage{530}
\lastpage{544}
\doi{10.3150/10-BEJ290}

\makeatletter

\newcommand{\mean}{\mathbb{E}}
\newcommand{\prob}{\mathbb{P}}
\newcommand{\var}{\operatorname{Var}}
\newcommand{\real}{\mathbb{R}}

\newcommand{\bn}{\mathbf{N}}
\newcommand{\bq}{\mathbf{Q}}
\newcommand{\bz}{\mathbf{Z}}
\newcommand{\bl}{\bolds\lambda}

\newcommand{\ca}{\mathcal{A}}
\newcommand{\cb}{\mathcal{B}}
\newcommand{\cf}{\mathcal{F}}
\newcommand{\ch}{\mathcal{H}}
\newcommand{\ci}{\mathcal{I}}
\newcommand{\cl}{\mathcal{L}}
\newcommand{\cu}{\mathcal{U}}

\newcommand{\G}{{\Gamma}}
\newcommand{\Pn}{\mathrm{Po}}

\newcommand{\Xii}{\Xi^{(i)}}

\newtheorem{Theorem}{Theorem}[section]
\newtheorem{Corollary}[Theorem]{Corollary}
\newtheorem{Lemma}[Theorem]{Lemma}
\newtheorem{Proposition}[Theorem]{Proposition}
\newremark{Remark}[Theorem]{Remark}
\newremark{Example}[Theorem]{Example}

\makeatother

\begin{document}
\begin{frontmatter}

\title{Poisson process approximation for dependent superposition of point processes}
\runtitle{Superposition of point processes}

\begin{aug}
\author[a]{\fnms{Louis H.Y.} \snm{Chen}\thanksref{a}\ead[label=e1]{matchyl@nus.edu.sg}}
\and
\author[b]{\fnms{Aihua} \snm{Xia}\thanksref{b}\ead[label=e2]{xia@ms.unimelb.edu.au}}

\runauthor{L.H.Y. Chen and A. Xia}

\address[a]{Department of Mathematics,
National University of Singapore,
10 Lower Kent Ridge Road,
Singapore 119076,
Republic of Singapore.
\printead{e1}}

\address[b]{Department of Mathematics and Statistics,
University of Melbourne,
VIC 3010, Australia.\\
\printead{e2}}
\end{aug}

\received{\smonth{1} \syear{2009}}
\revised{\smonth{4} \syear{2010}}

%
\begin{abstract}
Although the study of weak convergence of superpositions of point
processes to the Poisson process dates back to the work of Grigelionis
in 1963, it was only recently that Schuhmacher
[\textit{Stochastic Process. Appl.} \textbf{115} (2005) 1819--1837]
obtained error
bounds for the weak convergence. Schuhmacher considered dependent
{superposition}, truncated the individual point processes to 0--1 point
processes and then applied Stein's method to the latter. In this
paper, we adopt a different approach to the problem by using Palm
theory and Stein's method, thereby expressing the error bounds in terms
of the mean measures of the individual point processes, which is not
possible with Schuhmacher's approach. We consider locally dependent
{superposition} as a generalization of the locally dependent point
process introduced in Chen and Xia [\textit{Ann. Probab.} \textbf{32}
(2004) 2545--2569] and apply the main theorem
to the superposition of thinned point processes and of renewal
processes.
\end{abstract}

%

\end{frontmatter}

\section{Introduction}

The study of weak convergence of superpositions of point processes
dates back to Grigelionis \cite{Grigelionis63}
who proved that the {superposition} of independent sparse point
processes converges weakly to a Poisson process on the carrier space
$\real_+$. His result was subsequently extended to more general carrier
spaces by Goldman \cite{Goldman67} and Jagers \cite{Jagers72};
see \cite{Cinlar72}
and \cite{Brown78} for further discussion. It was further extended to
superpositions of dependent sparse point processes by Banys \cite
{Banis75,Banis85},
Kallenberg \cite{Kallenberg75}, Brown \cite{Brown79} and Banys \cite
{Banys80}. For a
systematic account of these developments, see \cite{Ka83}.

Surprisingly, it was only recently that error bounds for such
convergence of point processes were studied. Using Stein's method for
Poisson process approximation, as developed by Barbour \cite{Barbour88} and
Barbour and Brown \cite{BB92}, Schuhmacher \cite{Schuhmacher05a}
obtained an error bound
on the $d_2$ Wasserstein distance between a sum of weakly dependent
sparse point processes $\{\xi_{ni}, 1\le i\le k_n\}_{n\in\bn}$
and an approximating Poisson process. As he truncated the sparse point
processes to 0--1 point processes, as in the proof of Grigelionis'
theorem, his error bound contains the term $\sum_{i=1}^{k_n}\prob[\xi
_{ni}(B)\ge2]$, whose convergence to 0 for every bounded Borel subset
$B$ of the carrier space is a condition for Grigelionis' theorem to
hold. A consequence of such truncation is that the mean measure of the
approximating Poisson process is not equal to the sum of the mean
measures of the individual point processes.

In this paper, we adopt a different approach to Poisson process
approximation in which we do not use the truncation, but apply Palm
theory and express the error bounds in terms of the mean measures of
the individual sparse point processes. Such an approach also ensures
that the mean measure of the approximating Poisson process is equal to
the sum of the mean measures of the sparse point processes.

As in \cite{Schuhmacher05a}, we study the dependent superposition of
sparse point processes. But we consider only locally dependent
superposition, which is a natural extension of the point processes
$\sum I_i\delta_{U_i}$ studied in \cite{ChenXia04}, Section 4, where
$\delta_x$ is the point mass at $x$, the $U_i$'s are $\mathcal{S}$-valued
independent random elements with $\mathcal{S}$ a locally compact metric space,
the indicators $I_i$'s are locally dependent and the $I_i$'s are
independent of the $U_i$'s.

In our main theorem (Theorem \ref{maintheorem1}), with the help of
Brown, Weinberg and Xia \cite{BWX00}, Lemma~3.1, it is possible to
recover a
factor of order $1/\lambda$ from the term $1/(|\Xi^{(i)}|+1)$. Hence,
the error bound on the $d_2$ Wasserstein distance yields the so-called
Stein factor $1/\lambda$, by which approximation remains good for
large $\lambda
$, a feature always sought after for Poisson-type approximations. In
the error bound obtained by Schuhmacher \cite{Schuhmacher05a}, a
leading term does
not have the Stein factor; {see Remark \ref{rem4.4} for further details}.

Our main theorem and some corollaries are presented in Section
\ref{maintheorem}. Applications to thinned point processes and renewal
processes are given in Sections \ref{thinning} and \ref{renewal},
respectively.

\section{The main theorem}\label{maintheorem}

Throughout this paper, we assume that $\G$ is a
locally compact metric space with metric $d_0$ bounded by 1.
In estimating the error of Poisson process approximation to the
superposition of dependent point processes $\{\Xi_i, i\in\ci\}$ on the
carrier space $\G$ with $\ci$ a finite or countably infinite index set,
one natural approach is to partition the index set $\ci$ into $\{\{i\}
, \ci_i^s, \ci_i^w\}$, where $\ci_i^s$ is the set of indices of the
point processes which are strongly dependent on $\Xi_i$ and $\ci_i^w$
the set of the indices of the point processes which are weakly
dependent on $\Xi_i$; see \cite{Schuhmacher05a}. Another approach is to
divide the index set according to various levels of local dependence, a
successful structure for studying normal approximation; see \cite{ChenShao04}.
The latter approach has been generalized by Barbour and
Xia \cite{BarbourXia06} to randomly indexed sums with a particular
interest in
random variables resulting from integrating a random field
with respect to a point process.

Parallel to the local dependence structures defined in \cite{ChenShao04},
we introduce the following:

\begin{enumerate}[{[LD2]}]
\item[{[LD1]}] for each $i\in\ci$, there exists a neighborhood $A_i$
such that
$i\in A_i$ and $\Xi_i$ is independent of $\{\Xi_j, j\in A_i^c\}$;

\item[{[LD2]}] condition [LD1] holds and for each $i\in\ci$, there exists
a neighborhood $B_i$ such that $A_i\subset B_i$ and $\{\Xi_j, j\in
A_i\}$ is independent of $\{\Xi_j, i\in B_i^c\}$.
\end{enumerate}

The index set $\ci$ in [LD1] and [LD2] will be assumed to be finite or
countably infinite in this paper, although it may be as general as that
considered in \cite{BarbourXia06}. The superposition of $\{\Xi_i\dvtx
i\in\ci\}$ which satisfies the condition [LD1] is more general
than point processes of the form $\sum I_i\delta_{U_i}$, where the
$I_i$'s are locally dependent indicators with one level of dependent
neighborhoods in $\ci$ (i.e., the $I_i$'s satisfy [LD1] in
\cite{ChenShao04}, page 1986). Such a point process is a typical
example of locally dependent point processes defined in \cite{ChenShao04},
page 2548. Likewise, the superposition of $\{\Xi_i\dvtx i\in
\ci\}$ which satisfies the condition [LD2] is more general than point
processes of
the form $\sum I_i\delta_{U_i}$, where the $I_i$'s are locally
dependent indicators with two levels of
dependent neighborhoods in $\ci$ (i.e., the $I_i$'s satisfy [LD2] in
\cite{ChenShao04}, page 1986).

Three metrics will be used to describe the accuracy of Poisson process
approximation: the total variation metric for Poisson random variable
approximation $d_{\mathrm{tv}}$; the total variation metric for Poisson process
approximation $d_{\mathrm{TV}}$; and a Wasserstein metric $d_2$ (see
\cite{BHJ92}
or \cite{Xia05}).

To briefly define these metrics, let $\ch$ be the space of all finite
point process configurations on~$\G$, that is, each $\xi\in\ch$ is a
non-negative integer-valued finite measure on $\G$. Let $\mathcal{K}$
stand for the set of
$d_0$-Lipschitz functions $k\dvtx\G\rightarrow[-1,1]$ such that $\mid
k(\alpha)-k(\beta)\mid\le d_0(\alpha,\beta)$ for all
$\alpha,\beta\in\G.$ The first Wasserstein metric $d_1$ on $\ch$ is
defined by\looseness=-1
\[
d_1(\xi_1,\xi_2)=\cases{
0, &\quad if $|\xi_1|=|\xi_2|=0$,\cr
1, &\quad if $|\xi_1|\ne|\xi_2|$,\cr
\displaystyle|\xi_1|^{-1}\sup_{k\in\mathcal{K}} \biggl\vert\int k\,
\mathrm{d}\xi_1
-\int k\,\mathrm{d}\xi_2 \biggr\vert, &\quad if $|\xi_1|=|\xi_2|>0$,}
\]
where $|\xi_i|$ is the total mass of $\xi_i$.
A metric $d_1'$ equivalent to $d_1$ can be defined as follows (see
\cite{BX95}): for two configurations
$\xi_1=\sum_{i=1}^n\delta_{y_i}$ and $\xi_2=\sum_{i=1}^m\delta
_{z_i}$ with
$m\ge n$,
\[
d_1'(\xi_1,\xi_2)=\min_\pi\sum_{i=1}^nd_0\bigl(y_i,z_{\pi(i)}\bigr)+(m-n),
\]
where $\pi$ ranges over all permutations of $(1,\ldots,m)$. Both $d_1$
and $d_1'$ generate the weak topology on $\ch$ (see \cite{Xia05},
Proposition 4.2) and we use $\cb(\ch)$ to stand for the Borel $\sigma
$-algebra generated by the weak topology. Define three subsets of
real-valued functions on $\ch$:
$\cf_{\mathrm{tv}}=\{1_A(|\xi|)\dvtx A\subset\bz_+\}$, $\cf_{d_1}=\{f\dvtx  \vert
f(\xi
_1)-f(\xi_2)\vert\le d_1(\xi_1,\xi_2)\mbox{ for all }
\xi_1, \xi_2\in\ch\}$ and $\cf_{\mathrm{TV}}=\{1_A(\xi)\dvtx  A\in\cb(\ch
)\}$. The
pseudo-metric $d_{\mathrm{tv}}$ and the metrics $d_2$ and $d_{\mathrm
{TV}}$ are then
defined on probability measures on $\ch$ by
\begin{eqnarray*}
d_{\mathrm{tv}}(\bq_1,\bq_2)&=&\inf_{(X_1,X_2)}\prob(|X_1|\ne|X_2|)=\sup
_{f\in
\cf_{\mathrm{tv}}} \biggl|\int f\,\mathrm{d}\bq_1-\int f\,\mathrm{d}\bq_2 \biggr|,\\
d_2(\bq_1,\bq_2)&=&\inf_{(X_1,X_2)}\mean[d_1(X_1,X_2)]\\
&=&\sup_{f\in
\cf
_{d_1}} \biggl|\int f\,\mathrm{d}\bq_1-\int f\,\mathrm{d}\bq_2 \biggr|,\\
d_{\mathrm{TV}}(\bq_1,\bq_2)&=&\inf_{(X_1,X_2)}\prob(X_1\ne X_2)\\
&=&\sup
_{f\in
\cf
_{\mathrm{TV}}} \biggl|\int f\,\mathrm{d}\bq_1-\int f\,\mathrm{d}\bq_2 \biggr|,
\end{eqnarray*}
where the infima are taken over all
\textit{couplings} of $(X_1,X_2)$ such that $\cl(X_i)=\bq_i$, $i=1, 2$,
and the second equations are due to
the duality theorem; see \cite{Rachev91}, page 168.

To bound the error of Poisson process approximation, we {need} the Palm
distributions $Q_\alpha$ of a point process $X_2$ with respect to a point
process $X_1$ with finite mean measure $\nu$ at $\alpha$. When $X_1$
is a
\textit{simple} point process, that is, it has at most one point at each
location, the Palm distribution $Q_\alpha$ may be intuitively interpreted
as the conditional distribution of $X_2$ given that $X_1$ has a
point
at $\alpha$. More precisely, let $\cb(\G)$ denote the Borel $\sigma
$-algebra in $\G$ generated by the metric $d_0$ and define the Campbell
measure $C$ of $(X_1,X_2)$ on $\cb(\G)\times\cb(\ch)$:
\[
C(B\times M)=\mean[X_1(B)1_{X_2\in M}],\qquad B\in\cb(\G), M\in\cb
(\ch).
\]
Since the mean measure $\nu$ of $X_1$ is finite, by Kallenberg \cite{Ka83},
15.3.3, there exist probability measures $Q_\alpha$ on $\cb(\ch)$
such that
%
%
\begin{equation}\label{Palmdefinition1}
\mean[X_1(B)1_{X_2\in M} ]=\int_B
Q_\alpha(M)\nu
(\mathrm{d}\alpha)\qquad
\forall B\in\cb(\G), M\in\cb(\ch),
\end{equation}
which is equivalent to
\[
Q_\alpha(M)=\frac{\mean\{1_{[X_2\in M]}X_1(\mathrm{d}\alpha)\}}{\nu(\mathrm{d}\alpha
)}\qquad \forall M\in
\cb(\ch), \alpha\in\G\mbox{ $\nu$-a.s.;}
\]
%
%
see \cite{Ka83}, Section 10.1. It is possible to realize a family
of point processes $Y_\alpha$ on some probability space
such that $Y_\alpha\sim Q_\alpha$ and we say that $Y_\alpha$
is a \textit{Palm
process} of $X_2$ with respect to $X_1$ at $\alpha$. {Moreover, when
$X_1=X_2$, we call the point process $Y_\alpha-\delta_\alpha$ the
\textit
{reduced Palm process} of $X_2$ at $\alpha$; see \cite{Ka83},
Lemma 10.2.}

As noted in \cite{GX06}, when $\G$ is reduced to one point
only, the Palm distribution of $X_2$ (with respect to itself) is the
same as the size-biased distribution; general guidelines for the
construction of size-biased variables are investigated in \cite{GR96}.

\begin{Theorem}\label{maintheorem1} Let $\{\Xi_i,
i\in\ci\}$ be a collection of point processes on $\G$
with respective mean measures $\bl_i$, $i\in\ci$. Set
$\Xi=\sum_{i\in\ci}\Xi_i$ with mean measure denoted by $\bl$ and assume
that $\lambda:=\bl(\G)<\infty$. If \textup{[LD1]} holds, then
%
%
\begin{eqnarray}
\label{mainresult1}
d_{\mathrm{tv}}(\cl(\Xi),\Pn(\bl))&\le&\frac{1-\mathrm{e}^{-\lambda}}{\lambda
}\mean\sum
_{i\in\ci
}\int_\G\bigl\{ \bigl||V_i|-|V_{i,\alpha}| \bigr|+ \bigl||\Xi
_i|-|\Xi_{i,(\alpha)
}| \bigr| \bigr\}\bl_i(\mathrm{d}\alpha),\\
\label{mainresult2}
d_2(\cl(\Xi),\Pn(\bl))&\le&\mean\sum_{i\in\ci}
\biggl(\frac{3.5}{\lambda}+\frac{2.5}{|\Xii|+1} \biggr)\int_\G
d_1'(V_i,V_{i,\alpha})\bl_i(\mathrm{d}\alpha)\nonumber\\[-8pt]\\[-8pt]
&&{} +\sum_{i\in\ci}
\biggl(\frac{3.5}{\lambda}+\mean\frac{2.5}{|\Xii|+1} \biggr)\mean
\int
_\G
d_1'(\Xi_i,\Xi_{i,(\alpha)})\bl_i(\mathrm{d}\alpha),\nonumber\\
\label{mainresult3}
d_{\mathrm{TV}}(\cl(\Xi),\Pn(\bl))&\le&\mean\sum_{i\in\ci}
\int_\G\{\|V_i-V_{i,\alpha}\|+\|\Xi_i-\Xi_{i,(\alpha)}\|
\}\bl_i(\mathrm{d}\alpha),
\end{eqnarray}
where\vspace*{-2pt} $\Xii=\sum_{j\in A_i^c}\Xi_j$, $V_i=\sum_{j\in A_i\backslash
\{i\}
}\Xi_j$, $\Xi_{i,(\alpha)}$ is
the reduced Palm process of $\Xi_i$ at $\alpha$, $V_{i,\alpha}$ is
the Palm
process of $V_i$ with respect to $\Xi_i$ at $\alpha$ such that $\Xii
+V_{i,\alpha
}+\Xi_{i,(\alpha)}+\delta_\alpha$ is the Palm process of $\Xi$ with
respect to $\Xi
_i$ at $\alpha$ and $\|\cdot\|$ denotes the variation norm of signed measure.
Under the condition \textup{[LD2]}, (\ref{mainresult1}) and (\ref{mainresult3})
remain the same, but (\ref{mainresult2}) can be further reduced
to
%
%
\begin{eqnarray}\label{mainresult4}
d_2(\cl(\Xi),\Pn(\bl))&\le&\sum_{i\in\ci}
\biggl(\frac{3.5}{\lambda}+\mean\frac{2.5}{\sum_{j\in B_i^c}|\Xi
_j|+1}
\biggr)\mean\int_\G
d_1'(V_i,V_{i,\alpha})\bl_i(\mathrm{d}\alpha)\nonumber\\[-8pt]\\[-8pt]
&&{}+\sum_{i\in\ci}
\biggl(\frac{3.5}{\lambda}+\mean\frac{2.5}{|\Xii|+1} \biggr)\mean
\int
_\G
d_1'(\Xi_i,\Xi_{i,(\alpha)})\bl_i(\mathrm{d}\alpha)\nonumber\\
\label{xiaadd11}
&\le&\sum_{i\in\ci}
\biggl(\frac{3.5}{\lambda}+2.5\cdot\frac{\sqrt{\kappa_i(1+\kappa
_i/4)}+1+\kappa_i/2}{\sum_{j\in B_i^c}\lambda_j+1}
\biggr)\nonumber\\[-8pt]\\[-8pt]
&&{}\times\{\lambda_i\mean|V_i|
+\mean(|V_i|\cdot|\Xi_i|)+\lambda_i^2+\mean(|\Xi
_i|^2)-\lambda_i\},\nonumber
\end{eqnarray}
where $\lambda_i=\bl_i(\G)$ and
\[
\kappa_i=\frac{\sum_{j_1\in B_i^c}\sum_{j_2\in B_i^c\cap
A_{j_1}}\operatorname{cov}(|\Xi_{j_1}|,|\Xi_{j_2}|)}{\sum_{j\in
B_i^c}\lambda_j+1}.
\]
\end{Theorem}
%
%
\begin{pf}
We employ Stein's method for Poisson process
approximation, established in \cite{Barbour88} and \cite{BB92},
to prove the theorem. To this end, for a suitable
measurable function $h$ on $\ch$, let
\[
\ca h(\xi)=\int_\G[h(\xi+\delta_\alpha)-h(\xi)]\bl(\mathrm{d}\alpha
)+\int_\G
[h(\xi-\delta_x)-h(\xi)]\xi(\mathrm{d}x).
\]
Then $\ca$ defines a generator of the spatial immigration--death process
with immigration intensity $\bl$ and unit per capita death rate, and
the equilibrium distribution of the spatial immigration--death process
is $\Pn(\bl)$; see \cite{Xia05}, Section 3.2, for more details. The Stein
equation based on $\ca$ is
%
%
\begin{equation}\label{steinequation}
\ca h(\xi)=f(\xi)-\Pn(\bl)(f)
\end{equation}
with solution
\[
h_f(\xi)=-\int_0^\infty[\mean f(\bz_\xi(t))-\Pn(\bl)(f)]\,\mathrm{d}t,
\]
where $\{\bz_\xi(t), t\ge0\}$ is the spatial immigration--death
process with generator $\ca$ and initial configuration $\bz_\xi
(0)=\xi
$. To obtain bounds on the errors in the approximation, we need to define
\begin{eqnarray*}
\Delta h_f(\xi;x)&:=&h_f(\xi+\delta_x)-h_f(\xi),\\
\Delta^2h_f(\xi;x,y)&:=&\Delta h_f(\xi+\delta_x;y)-\Delta h_f(\xi
;y),\\
\Delta^2h_f(\xi,\eta;x)&:=&\Delta h_f(\xi;x)-\Delta h_f(\eta;x),
\end{eqnarray*}
for corresponding test functions $f$. Xia \cite{Xia05}, Propositions
5.6 and
5.12 (see \cite{BE83,BB92})
and Lemma~5.26, state that, for all $x,y\in\G$,
%
%
\begin{eqnarray}
\label{steinconstanttv}
|\Delta^2h_f(\xi;x,y)|&\le&\frac{1-\mathrm{e}^{-\lambda}}{\lambda}\qquad
\forall
f\in\cf_{\mathrm{tv}},\\
\label{steinconstantprocesstv}
|\Delta^2h_f(\xi;x,y)|&\le&1\qquad \forall f\in\cf_{\mathrm{TV}},\\
\label{steinconstantd1}
|\Delta^2h_f(\xi,\eta;x)|&\le&\biggl(\frac{3.5}{\lambda}+\frac
{2.5}{|\eta
|\wedge|\xi|+1} \biggr)d_1'(\xi,\eta)\qquad \forall f\in\cf
_{d_1}.
\end{eqnarray}
%

Now, since $\Xii+V_{i,\alpha}+\Xi_{i,(\alpha)}+\delta_\alpha$ is the
Palm process of $\Xi
$ with respect to $\Xi_i$ at $\alpha$, it follows from (\ref
{Palmdefinition1}) that
\begin{eqnarray*}
\mean
\int_\G[h(\Xi)-h(\Xi-\delta_\alpha)]\Xi(\mathrm{d}\alpha)
&=&\sum_{i\in\ci}\mean\int_\G[h(\Xi)-h(\Xi-\delta_\alpha)]\Xi
_i(\mathrm{d}\alpha)\\
&=&\sum_{i\in\ci}\mean\int_\G\Delta h\bigl(\Xii+V_{i,\alpha}+\Xi
_{i,(\alpha)};\alpha\bigr)\bl
_i(\mathrm{d}\alpha).
\end{eqnarray*}
On the other hand, by the Stein equation {(\ref{steinequation})}, we have
%
%
\begin{eqnarray}
&& |\mean f(\Xi)-\Pn(\bl)(f) |\nonumber\\
&&\quad= \biggl|\mean
\int_\G[h_f(\Xi+\delta_\alpha)-h_f(\Xi)]\bl(\mathrm{d}\alpha)+\mean
\int_\G[h_f(\Xi-\delta_x)-h_f(\Xi)]\Xi(\mathrm{d}x) \biggr|\nonumber\\
\label{mainresultproof1}
&&\quad= \biggl|\sum_{i\in\ci}\mean\int_\G\bigl\{\Delta h_f(\Xi;\alpha
)-\Delta h_f\bigl(\Xii
+V_{i,\alpha}+\Xi_{i,(\alpha)};\alpha\bigr)\bigr\}\bl_i(\mathrm{d}\alpha) \biggr|
\\
\label{mainresultproof2}
&&\quad\le\sum_{i\in\ci}\mean\int_\G\bigl\{\bigl|\Delta h_f\bigl(\Xii+V_i+\Xi
_i;\alpha\bigr)-\Delta
h_f\bigl(\Xii+V_{i,\alpha}+\Xi_i;\alpha\bigr)\bigr|\nonumber\\[-8pt]\\[-8pt]
&&\qquad\hspace*{36pt}{}+\bigl|\Delta h_f\bigl(\Xii+V_{i,\alpha}+\Xi_i;\alpha\bigr)-\Delta h_f\bigl(\Xii
+V_{i,\alpha}+\Xi
_{i,(\alpha)};\alpha\bigr)\bigr|\bigr\}\bl_i(\mathrm{d}\alpha).\qquad\quad\nonumber
\end{eqnarray}

To prove (\ref{mainresult1}), we note that the test functions $f\in
\cf
_{\mathrm{tv}}$ satisfy $f(\xi)=f(|\xi|\delta_z)$ for a fixed point
$z\in\G
$ and
so we have $h_f(\xi)=h_f(|\xi|\delta_z)$. Hence, for all $\eta,
\xi
_1, \xi_2\in\ch$,
%
%
\begin{eqnarray}\label{mainresultproof3}
&&|\Delta h_f(\eta+\xi_1;\alpha)-\Delta h_f(\eta+\xi_2;\alpha
)|\nonumber\\
&&\quad=\bigl|\Delta h_f\bigl(\eta+(|\xi_1|\vee|\xi_2|)\delta_z;\alpha\bigr)-\Delta
h_f\bigl(\eta
+(|\xi_1|\wedge|\xi_2|)\delta_z;\alpha\bigr)\bigr|\\
&&\quad\le\sum_{j=1}^{ ||\xi_1|- |\xi_2| |} \bigl|\Delta^2
h_f\bigl(\eta+(|\xi
_1|\wedge|\xi_2|+j-1)\delta_z;z,\alpha\bigr)\bigr|\le\bigl||\xi_1|- |\xi_2|
\bigr|\frac{1-\mathrm{e}^{-\lambda}}{\lambda},\nonumber
\end{eqnarray}
where the last inequality is due to (\ref{steinconstanttv}). Combining
(\ref{mainresultproof3}) with (\ref{mainresultproof2}) yields (\ref
{mainresult1}).

Next, (\ref{steinconstantd1}) and (\ref{mainresultproof1}) imply
that for
$f\in\cf_{d_1}$,
\[
|\mean f(\Xi)-\Pn(\bl)(f) | \le\sum_{i\in\ci}\mean\int_\G
\biggl(\frac{3.5}{\lambda}+\frac{2.5}{|\Xii|+1} \biggr)d_1'(V_i+\Xi
_i,V_{i,\alpha
}+\Xi_{i,(\alpha)})\bl_i(\mathrm{d}\alpha).
\]
Because $d_1'(V_i+\Xi_i,V_{i,\alpha}+\Xi_{i,(\alpha)})\le
d_1'(V_i,V_{i,\alpha
})+d_1'(\Xi_i,\Xi_{i,(\alpha)})$ and, for each $i\in\ci$, $\Xi_i$ is
independent of $\Xii$, (\ref{mainresult2}) follows. On the other hand,
due to the independence between $\{V_i,\Xi_i\}$ and $\{\Xi_j, j\in
B_i^c\}$ {implied by [LD2]}, (\ref{mainresult4}) is immediate. To prove
(\ref{xiaadd11}), one can verify that
\[
\var\biggl(\sum_{j\in B_i^c}|\Xi_j| \biggr)=\sum_{j_1\in B_i^c}\sum
_{j_2\in B_i^c\cap A_{j_1}}\operatorname{cov}(|\Xi_{j_1}|,|\Xi_{j_2}|)
\]
and that $\mean\sum_{j\in B_i^c}|\Xi_j|=\sum_{j\in B_i^c}\lambda
_j$. Hence,
(\ref{xiaadd11}) follows from Lemma 3.1 in \cite{BWX00}
and the facts that
$d_1'(V_i,V_{i,\alpha})\le|V_i|+|V_{i,\alpha}|$ and $d_1'(\Xi_i,\Xi
_{i,(\alpha)})\le
|\Xi_i|+|\Xi_{i,(\alpha)}|$.

Finally, we show (\ref{mainresult3}). For $\xi_1, \xi_2\in\ch$,
we define
\[
\xi_1\wedge\xi_2=\sum_{j=1}^k(a_{1j}\wedge a_{2j})\delta_{x_j},
\]
where $\{x_1,\dots,x_k\}$ is the support of the point measure $\xi
_1+\xi
_2$, so that $\xi_i=\sum_{j=1}^ka_{ij}\delta_{x_j}$ for $i=1,2$ with
the $a_{ij}$'s being non-negative integers. Then, for all $f\in\cf
_{\mathrm{TV}}, \eta, \xi_1, \xi_2\in\ch$,
%
%
\begin{eqnarray}\label{mainresultproof4}
&&|\Delta h_f(\eta+\xi_1;\alpha)-\Delta h_f(\eta+\xi_2;\alpha
)|\nonumber\\
&&\quad\le|\Delta h_f(\eta+\xi_1;\alpha)-\Delta h_f(\eta+\xi_1\wedge
\xi_2;\alpha
)|\nonumber\\
&&\qquad{}+|\Delta h_f(\eta+\xi_2;\alpha)-\Delta h_f(\eta+\xi_1\wedge\xi
_2;\alpha)|\\
&&\quad\le(|\xi_1|-|\xi_1\wedge\xi_2|)+(|\xi_2|-|\xi_1\wedge\xi
_2|)\nonumber\\
&&\quad=\|\xi_1-\xi_2\|,\nonumber
\end{eqnarray}
where the last inequality is due to (\ref{steinconstantprocesstv}).
Applying (\ref{mainresultproof4}) in
(\ref{mainresultproof2}), we obtain (\ref{mainresult3}).
\end{pf}
\begin{Corollary} \label{Corollaryindeptcase} With the notation of
Theorem \ref{maintheorem1},
if $\{\Xi_i, i\in\ci\}$ are all independent, then
%
%
\begin{eqnarray}
\label{mainresult6}
d_{\mathrm{tv}}(\cl(\Xi),\Pn(\bl))&\le&\frac{1-\mathrm{e}^{-\lambda}}{\lambda
}\mean\sum
_{i\in\ci
}\int_\G\bigl||\Xi_i|-|\Xi_{i,(\alpha)}| \bigr|\bl_i(\mathrm{d}\alpha
),\\
\label{xiaadd12}
d_2(\cl(\Xi),\Pn(\bl))&\le&\sum_{i\in\ci} \biggl(\frac
{3.5}{\lambda
}+\mean\frac
{2.5}{\sum_{j\ne i}|\Xi_j|+1} \biggr)\mean\int_\G d_1'(\Xi_i,\Xi
_{i,(\alpha)
})\bl_i(\mathrm{d}\alpha)\\
\label{mainresult5}
&\le& \biggl(\frac{3.5}{\lambda}+2.5\cdot\frac{\sqrt{\kappa
(1+\kappa
/4)}+1+\kappa/2}{\lambda-\max_{j\in\ci}\lambda_j+1}
\biggr) \sum
_{i\in\ci
}\{\lambda
_i^2+\mean(|\Xi_i|^2 )-\lambda_i\},\hspace*{33pt}\\
\label{mainresult7}
d_{\mathrm{TV}}(\cl(\Xi),\Pn(\bl))&\le&\mean\sum_{i\in\ci}
\int_\G\|\Xi_i-\Xi_{i,(\alpha)}\|\bl_i(\mathrm{d}\alpha),
\end{eqnarray}
where $\kappa=\frac{\sum_{i\in\ci}\var(|\Xi_i|)}{\lambda-\max
_{j\in
\ci}\lambda_j+1}$.
\end{Corollary}
\begin{pf}
Let $A_i=B_i=\{i\}$, then (\ref{mainresult6})--(\ref{mainresult7}) follow
from (\ref{mainresult1}), (\ref{mainresult4}), (\ref{xiaadd11}) and
(\ref{mainresult3}),
respectively.
\end{pf}
\begin{Corollary}[(cf. \cite{ChenXia04}, Theorem 4.1)] \label
{corollarychenxia04}
Let $\{I_i, i\in\ci\}$ be dependent
indicators with $\ci$ a finite or countably infinite index set and let
$\{\cu_i, i\in\ci\}$ be $\G$-valued independent random elements
independent of $\{I_i, i\in\ci\}$. Define $\Xi=\sum_{i\in\ci
}I_i\delta_{\cu
_i}$ with mean measure $\bl$, let $\mean I_i=p_i$ and assume that
$\lambda
=\sum_{i\in\ci}p_i<\infty$. For each $i\in\ci$, let $A_i$ be the
set of
indices of those $I_j$'s which are dependent on $I_i$, that is, $I_i$
is independent of $\{I_j\dvtx  j\in A_i^c\}$. Then,
%
%
\begin{eqnarray}\label{xiaadd1}
d_{\mathrm{tv}}(\cl(\Xi), \Pn(\bl))
&\le&
\frac{1-\mathrm{e}^{-\lambda}}{\lambda}\sum_{i\in\ci} \biggl\{\sum_{j\in
A_i\backslash\{i\}
}\mean I_iI_j+\sum_{j\in
A_i}p_ip_j \biggr\},\\
\label{xiaadd2}
d_2(\cl(\Xi), \Pn(\bl))
&\le&
\mean\sum_{i\in\ci}\sum_{j\in A_i\backslash\{i\}} \biggl(\frac
{3.5}{\lambda
}+\frac{2.5}{{S}_i+1} \biggr)I_iI_j\nonumber\\[-8pt]\\[-8pt]
&&{}+\sum_{i\in\ci}\sum_{j\in
A_i} \biggl(\frac{3.5}{\lambda}+\mean\biggl[\frac
{2.5}{{S}_i+1}
\Big|I_j=1 \biggr] \biggr)p_ip_j,\nonumber\\
\label{xiaadd3}
d_{\mathrm{TV}}(\cl(\Xi), \Pn(\bl))
&\le&\sum_{i\in\ci} \biggl\{\sum_{j\in A_i\backslash\{i\}}\mean
I_iI_j+\sum_{j\in
A_i}p_ip_j \biggr\},
\end{eqnarray}
where ${S}_i=\sum_{j\notin A_i}I_j$. For each $i\in\ci$, let $B_i$ be
the set of indices of those $I_l$'s which are dependent on $\{I_j,
j\in A_i\}$ so that $\{I_j\dvtx  j\in A_i\}$ is independent of $\{I_l\dvtx
l\in B_i^c\}$. Then, (\ref{xiaadd1}) and (\ref{xiaadd3}) remain the same,
but (\ref{xiaadd2}) can be further reduced to
%
%
\begin{eqnarray}\label{xiaadd4}
d_2(\cl(\Xi), \Pn(\bl))
&\le&
\sum_{i\in\ci}\sum_{j\in A_i\backslash\{i\}} \biggl(\frac
{3.5}{\lambda
}+\mean
\frac{2.5}{W_i+1} \biggr)\mean(I_iI_j)\nonumber\\[-8pt]\\[-8pt]
&&{} +\sum_{i\in\ci}\sum_{j\in
A_i} \biggl(\frac{3.5}{\lambda}+\mean\frac{2.5}{W_i+1}
\biggr)p_ip_j\nonumber\\
\label{xiaadd5}
&\le&\sum_{i\in\ci} \biggl(\frac{3.5}{\lambda}+2.5\cdot\frac
{\sqrt
{\kappa
_i(1+\kappa_i/4)}+1+\kappa_i/2}{\sum_{j\in B_i^c}p_j+1}
\biggr)\nonumber\\[-8pt]\\[-8pt]
&&{} \times
\biggl(\sum_{j\in A_i\backslash\{i\}}\mean I_iI_j+\sum_{j\in
A_i}p_ip_j \biggr),\nonumber
\end{eqnarray}
where $W_i=\sum_{j\notin B_i}I_j$ and
\[
\kappa_i=\frac{\sum_{j_1\in B_i^c}\sum_{j_2\in B_i^c\cap
A_{j_1}}\operatorname{cov}(I_{j_1},I_{j_2})}{\sum_{j\in B_i^c}p_j+1}.
\]
\end{Corollary}
\begin{pf}
If we set $\Xi_i=I_i\delta_{\cu_i}$, $i\in
\ci$, then
$\Xi_i$ is independent of $\{\Xi_j\dvtx  j\notin A_i\}$, so [LD1] holds,
$\Xi_{i,(\alpha)}=0$, and the claims (\ref{xiaadd1})--(\ref{xiaadd3})
follow from (\ref{mainresult1})--(\ref{mainresult3}),
respectively. On the other hand, $\{\Xi_j\dvtx
j\in A_i\}$
is independent of
$\{\Xi_j\dvtx  j\notin B_i\}$, so [LD2] holds and (\ref{xiaadd4}) and
(\ref{xiaadd5}) are direct consequences of (\ref{mainresult4}) and
(\ref{xiaadd11}).
\end{pf}

A typical example of Poisson process approximation is that of
the Bernoulli process defined as follows (see \cite{Xia05}, Section
6.1,
for further discussion).
Let $I_1, \dots,I_n$ be independent indicators with
\[
\prob(I_i=1)=1-\prob(I_i=0)=p_i,\qquad i=1,\dots,n.
\]
Let $\Gamma=[0,1]$, $\Xi=\sum_{i=1}^nI_i\delta_{i/n}$ and $\bl
=\sum
_{i=1}^np_i\delta_{i/n}$ be the
mean measure of $\Xi$. {If we set} $\Xi_i=I_i\delta_{i/n}$,
$i=1,\dots
,n$, {then the reduced Palm process of $\Xi_i$ at $\alpha\in\G$ is
$\Xi
_{i,(\alpha)}=0$ and the Palm distribution of $\Xi_j$ with respect to point
process $\Xi_i$ at $\alpha$ for $j\ne i$ is the same as that of~$\Xi_j$.
Hence,} Corollary \ref{Corollaryindeptcase}, together with
(\ref{xiaadd12}) and \cite{ChenXia04}, Proposition 4.5, can be used to
obtain immediately the following {(known)} result.

\begin{Example}[(\cite{Xia05}, Section 6.1)]
For the Bernoulli process $\Xi$ on $\Gamma=[0,1]$ with mean measure
$\bl$,
\begin{eqnarray*}
d_{\mathrm{tv}}(\cl(\Xi),\Pn(\bl))&\le&\frac{1-\mathrm{e}^{-\lambda}}{\lambda
}\sum
_{i=1}^np_i^2,\\
d_{\mathrm{TV}}(\cl(\Xi),\Pn(\bl))&\le&\sum_{i=1}^np_i^2,\\
d_2(\cl(\Xi),\Pn(\bl))&\le&\frac{6}{\lambda-\max_{1\le i\le
n}p_i}\sum_{i=1}^np_i^2.
\end{eqnarray*}
\end{Example}

\begin{Example}
Throw $n$ points uniformly
and independently onto the interval $[0,n]$ and let $\Xi$ be the
configuration of the points on $[0,T]:=\G$ with $n\gg T$ and $\bl$ be
the mean measure of $\Xi$. Then,
\[
d_2(\cl(\Xi),\Pn(\bl))\le\frac{6T}{n-1}.
\]
\end{Example}
\begin{pf}
Let $I_i=1$ if the $i$th point is in $\G$ and
$0$ if it is not in $\G$. The configuration of the $i$th point on $\G$
can then be written as
$\Xi_i=I_i\delta_{\cu_i}$ and $\Xi=\sum_{i=1}^n\Xi_i$, where the
$\cu_i$'s
are independent and identically distributed uniform random variables on
$\G$ and are independent of the $I_i$'s. Noting that the reduced Palm
process $\Xi_{i,(\alpha)}=0$, we obtain the bound by applying {(\ref{xiaadd4})}
with $p_i=\prob(I_i=1)=T/n$ and \cite{ChenXia04}, Proposition 4.5.
\end{pf}

\section{Superposition of thinned dependent point processes}\label{thinning}

Assume that $q$ is a measurable retention function on $\G$ and $X$ is a
point process on $\G$. For a realization $X(\omega)$ of $X$, we thin
its points as follows. For each point of $X(\omega)$ at $\alpha$, it is
retained with probability $q(\alpha)$ and discarded with probability
$1-q(\alpha
)$, independently of the other points; see \cite{DaVerJon88},
page 554, for dependent thinning and \cite{Schuhmacher05b} for
discussions of more thinning strategies. The thinned configuration is
denoted by $X_q(\omega)$. For retention functions $q_1,q_2,\dots,q_n$,
let $\sum_{i=1}^nX'_{q_i}$ be the process arising from the
superposition of independent realizations of $X_{q_1},X_{q_2},\dots
,X_{q_n}$, that is, $X'_{q_1},X'_{q_2},\dots,X'_{q_n}$ are independent
and $\cl(X'_{q_i} )=\cl(X_{q_i} )$ for
$i=1,\dots
,n$. Fichtner \cite{Fichtner77} showed that a sequence of such superpositions,
obtained from the rows of an infinitesimal array of retention
functions, converges to a Poisson process under standard conditions;
see also \cite{Ka83}, Exercise 8.8).
Serfozo \cite{Serfozo84} presented convergence theorems for sums of dependent
point processes that are randomly thinned by a two-step procedure which
deletes each entire point process with a given probability and for each
retained point process, points are deleted or retained according to
another thinning strategy.
Necessary and sufficient conditions are given for a sum of two-step
thinned point processes to converge in distribution and the limit is
shown to be a Cox process; see also \cite{Fichtner75} and \cite{Liese80}.

For simplicity, we assume that $\{\Xi_i, i\in\ci\}$ is a locally
dependent collection of point processes (satisfying [LD1]) on a locally
compact metric space $\G$ with metric $d_0$ bounded by 1. For each
point of $\Xi_i$, we delete the point with probability $1-p$ and retain
it with probability $p$, independent of the others. The thinned point
process is denoted by $\Xi_i^p$, $i\in\ci$, and, in general, for each
point process $X$, we use $X^p$ to denote its thinned process. Let $\Xi
^p=\sum_{i\in\ci}\Xi_i^p$. As before, we define $A_i$ to be the
collection of indices $j$ of the point processes $\Xi_j$ which are
dependent on $\Xi_i$, that is, $\Xi_i$ is independent of $\{\Xi_j,
j\in A_i^c\}$.

\begin{Theorem}\label{mainthinningsresult} Let $\bolds\mu_i$ be the mean
measure of $\Xi_i$, $\mu_i=\bolds\mu_i(\G)=\mean(|\Xi_i|)$,
$i\in\ci$, and
assume that $\lambda=\sum_{i\in\ci}\mu_i<\infty$. The mean
measure of
$\Xi
^p$ is then $\bl^p=p\sum_{i\in\ci}\bolds\mu_i$ and
%
%
\begin{eqnarray}\label{thinningresult1}
d_{\mathrm{tv}}(\cl(\Xi^p),\Pn(\bl^p))&\le& p \biggl(1\wedge\frac
{1}{\lambda
}
\biggr)\mean\sum_{i\in\ci}
\{[|V_i|+|\Xi_i|]\lambda_i+[|V_i|+|\Xi_i|-1]|\Xi_i|\},\\
\label{thinningresult2}
d_2(\cl(\Xi^p),\Pn(\bl^p))&\le& p\mean\sum_{i\in\ci}
\biggl(\frac{3.5}{\lambda}+\frac{2.5}{ |\sum_{j\in A_i^c}\Xi
_j
|+1} \biggr)\nonumber\\[-8pt]\\[-8pt]
&&\hspace*{27pt}{}\times\{[|V_i|+|\Xi_i|]\lambda_i+[|V_i|+|\Xi_i|-1]|\Xi_i|\},
\nonumber\\
\label{thinningresult3}
d_{\mathrm{TV}}(\cl(\Xi^p),\Pn(\bl^p))&\le& p\mean\sum_{i\in\ci}
\{[|V_i|+|\Xi_i|]\lambda_i+[|V_i|+|\Xi_i|-1]|\Xi_i|\}.
\end{eqnarray}
\end{Theorem}
\begin{pf}
We prove only (\ref{thinningresult2}), as the
proofs of (\ref{thinningresult1}) and (\ref{thinningresult3}) are similar
to that of (\ref{thinningresult2}). By conditioning on the
configurations, we have, for each Borel set $B\subset\G$,
\[
\mean[\Xi_i^p(B)]=\mean\{\mean[\Xi_i^p(B)|\Xi_i]\}=\mean[\Xi_i(B)p],
\]
which implies that the mean measure of $\Xi_i^p$ is $\bl_i^p=p\bolds
\mu_i$
and hence that $\bl^p=p\sum_{i\in\ci}\bolds\mu_i$. By (\ref
{mainresult2}) and
the fact that $d_1'(\xi_1,\xi_2)\le|\xi_1|+|\xi_2|$, we obtain
\begin{eqnarray*}
&&d_2(\cl(\Xi^p),\Pn(\bl^p))\\
&&\quad\le
\mean\sum_{i\in\ci}
\biggl(\frac{3.5}{p\lambda}+\frac{2.5}{ |\sum_{j\in A_i^c}\Xi
_j^p
|+1} \biggr)\int_\G[|V_i^p|+|V_{i,\alpha}^p|+|\Xi_i^p|+|\Xi
_{i,(\alpha)}^p|]\bl
_i^p(\mathrm{d}\alpha)\\
&&\quad\le
\mean\sum_{i\in\ci}
\biggl(\frac{3.5}{p\lambda}+\frac{2.5}{ |\sum_{j\in A_i^c}\Xi
_j^p
|+1} \biggr)\int_\G\{[|V_i^p|+|\Xi_i^p|]\bl_i^p(\mathrm{d}\alpha
)+[|V_i^p|+|\Xi
_i^p|-1]\Xi_i^p(\mathrm{d}\alpha)\}\\
&&\quad\le
\mean\sum_{i\in\ci}
\biggl(\frac{3.5}{p\lambda}+\frac{2.5}{ |\sum_{j\in A_i^c}\Xi
_j^p
|+1} \biggr)\{[|V_i^p|+|\Xi_i^p|]\lambda_ip+[|V_i^p|+|\Xi
_i^p|-1]|\Xi
_i^p|\}.
\end{eqnarray*}
Since the points are thinned independently, we can condition on the
configuration of $\{\Xi_i, i\in\ci\}$. Noting that for
$Z\sim\operatorname{Binomial}(n,p)$, $\mean\frac{1}{Z+1}\le\frac
{1}{(n+1)p}$ and $\mean
[(X-1)X]=n(n-1)p^2$, we obtain
\begin{eqnarray*}
&&d_2(\cl(\Xi^p),\Pn(\bl^p))\\
&&\quad\le\mean\sum_{i\in\ci}
\biggl(\frac{3.5}{p\lambda}+\frac{2.5}{p ( |\sum_{j\in
A_i^c}\Xi
_j |+1 )} \biggr)\{[|V_i|+|\Xi_i|]\lambda_i+[|V_i|+|\Xi
_i|-1]|\Xi
_i|\}p^2.
\end{eqnarray*}
This completes the proof of (\ref{thinningresult2}).
\end{pf}

\begin{Remark}
Serfozo \cite{Serfozo84}, Example 3.6, obtained the rate $p$
for the convergence of a sum of thinned point processes to a Poisson
process. Theorem \ref{mainthinningsresult} shows that the rate $p$ is
valid for all of the three metrics used.
\end{Remark}

\section{Superposition of renewal processes}\label{renewal}

Viswanathan \cite{Viswanathan92}, page 290, states that if $\{\Xi_i,
1\le i\le n\}$
are independent renewal processes on $[0,T]$, each representing the
process of calls generated by a subscriber, then the total number of
calls can be modeled by a Poisson process. In this section, we quantify
this statement by giving an error bound for Poisson process
approximation to the sum of independent sparse renewal processes. We
begin with a technical lemma.

\begin{Lemma}
Let $\eta\sim G,$ $\xi_i\sim F, i\ge1$, be independent non-negative
random variables and define
\[
N_t=\max\{n\dvtx  \eta+\xi_1+\cdots+\xi_{n-1}\le t\},\qquad t\ge0.
\]
Then,
%
%
\begin{eqnarray}\label{renewallemma1}
G(t) &\le& \mean(N_t)\le\frac{G(t)}{1-F(t)},\\
\label{renewallemma2}
\mean(N_t^2)-\mean(N_t)&\le&\frac{2F(t)\mean(N_t)}{1-F(t)}\le
\frac
{2F(t)G(t)}{(1-F(t))^2}.
\end{eqnarray}
\end{Lemma}
\begin{pf}
Let $V(t)=\mean(N_t)$. The renewal equation gives
%
%
\begin{equation}\label{renewal1}
V(t)=G(t)+\int_0^t V(t-s)\,\mathrm{d}F(s)
\le G(t)+V(t)F(t),
\end{equation}
which implies (\ref{renewallemma1}). For (\ref{renewallemma2}), define
$V_2(t)=\mean[N_t(N_t+1)]$. Then, using the same arguments as for
proving the renewal equation,
\[
V_2(t)=2V(t)+\int_0^tV_2(t-s)\,\mathrm{d}F(s).
\]
This implies that $V_2(t)\le2V(t)+V_2(t)F(t)$, which, in turn, implies
that
\[
V_2(t)\le\frac{2V(t)}{1-F(t)}.
\]
Since
\[
\mean(N_t^2)-\mean(N_t)=V_2(t)-2V(t)=\int_0^tV_2(t-s)\,\mathrm{d}F(s)\le
V_2(t)F(t)\le\frac{2F(t)V(t)}{1-F(t)},
\]
(\ref{renewallemma2}) follows from (\ref{renewallemma1}).
\end{pf}
\begin{Theorem}\label{renewaltheorem} Suppose that $\{\Xi_i, 1\le
i\le
n\}$ are independent renewal processes on $[0,T]$ with the first
arrival time of $\Xi_i$ having distribution $G_i$ and its inter-arrival
time having distribution $F_i$. Let $\Xi=\sum_{i=1}^n\Xi_i$ and $\bl$
be its mean measure. Then,
%
%
\begin{equation}\label{renewalmainresult}
d_2(\cl(\Xi),\Pn(\bl))\le\frac{6\sum
_{i=1}^n[2F_i(T)+G_i(T)]G_i(T)}{ (\sum_{i=1}^n G_i(T)-\max
_jG_j(T) )(1-F_i(T))^2}.
\end{equation}
\end{Theorem}
\begin{pf}
We view a renewal process as a point process
whose points occur at the \mbox{renewal} times. For a renewal process $X$ with
renewal times $\tau_1\le\tau_2\le\cdots$, we further define $X'=\delta_{\tau_1}$.
Since $\lambda_i=\mean(|\Xi_i|)$, it follows from (\ref{xiaadd12}) that
%
%
\begin{eqnarray}\label{xiaadd21}
d_2(\cl(\Xi),\Pn(\bl))&\le&\sum_{i=1}^n \biggl(\frac{3.5}{\lambda
}+\mean\frac
{2.5}{\sum_{j\ne i}|\Xi_j|+1} \biggr)[\lambda_i^2+\mean(|\Xi
_i|^2)-\lambda
_i]\nonumber\\[-8pt]\\[-8pt]
&\le&\sum_{i=1}^n \biggl(\frac{3.5}{\lambda}+\mean\frac{2.5}{\sum
_{j\ne i}|\Xi
_j'|+1} \biggr) [\lambda_i^2+\mean(|\Xi_i|^2)-\lambda_i].\nonumber
\end{eqnarray}
However, applying Proposition 4.5 of \cite{ChenShao04} gives
\[
\mean\frac{1}{\sum_{j\ne i}|\Xi_j'|+1}\le\frac{1}{\sum_{j\ne
i}\mean
|\Xi_j'|}
=\frac{1}{\sum_{j\ne i}G_j(T)}
\]
and using (\ref{renewallemma1}), we obtain
\[
\lambda\ge\sum_{i=1}^n G_i(T).
\]
By combining (\ref{xiaadd21}), (\ref{renewallemma1}) and (\ref
{renewallemma2}), we obtain (\ref{renewalmainresult}).
\end{pf}
\begin{Remark}\label{remark4.3}
If $\{\Xi_i, 1\le i\le n\}$ are
independent and identically distributed stationary renewal processes on
$[0,T]$ with the successive inter-arrival time distribution $F$, then
\[
d_2(\cl(\Xi),\Pn(\bl))\le\frac{6n[2F(T)+G(T)]}{(n-1)(1-F(T))^2},
\]
where $G(t)=\int_0^t(1-F(s))\,\mathrm{d}s/\int_0^\infty(1-F(s))\,\mathrm{d}s$; see \cite
{DaVerJon88}, page 71.
\end{Remark}

\begin{Remark}\label{rem4.4}
An application of \cite{Schuhmacher05a}, Theorem 2.1,
to the sum of the renewal processes $\{\Xi_i, 1\le i\le n\}$ in
Remark \ref{remark4.3} with the natural partition $\{\{i\},\varnothing
,\{
1,\dots,i-1,i+1,\dots,n\}\}$ for each $1\le i\le n$ will give an
error bound
\[
n[F(T)+G(T)]+\theta G(T)(1+\ln^+n),
\]
where $\theta$ is a constant. The first term of the bound increases
linearly in $n$ and the bound is clearly not as sharp as the bound in
Remark \ref{remark4.3}.
\end{Remark}

Since the thinned process $X^p$ of a renewal process $X$ with mean
measure $\bolds\mu$ is still a renewal process (see \cite{DaVerJon88},
pages 75--76) with mean measure $\bolds\mu
^p=p\bolds\mu$ (see
the proof of Theorem \ref{mainthinningsresult}), {a~repetition of the
proof of Theorem \ref{renewaltheorem} yields the following proposition}.

\begin{Proposition} Suppose that $\{\Xi_i, 1\le i\le n\}$ are
independent renewal processes on $[0,T]$ with the first arrival time of
$\Xi_i$ having distribution $G_i$ and its inter-arrival time having
distribution $F_i$. Let $\Xi_i^p$ be the thinned point process obtained
from $\Xi_i$ by deleting each point with probability $1-p$ and
retaining it with probability $p$, independently of the other points.
Let $\Xi^p=\sum_{i=1}^n\Xi_i^p$ and $\bl^p$ be its mean measure. Then,
\[
d_2(\cl(\Xi^p),\Pn(\bl^p))\le\frac{6p\sum
_{i=1}^n[2F_i(T)+G_i(T)]G_i(T)}{ (\sum_{i=1}^n G_i(T)-\max
_jG_j(T) )(1-F_i(T))^2}.
\]
\end{Proposition}

\section*{Acknowledgements}
The research of the L.H.Y. Chen was partially supported by Grant
C-389-000-010-101 at the National University of Singapore and the Belz
Fund of the University of Melbourne. The research of A. Xia was
partially supported by the ARC Centre of Excellence for Mathematics and
Statistics of Complex Systems.


%
\printhistory


\begin{thebibliography}{00}

\bibitem{Banis75}
Banys, R. (1975).
The convergence of sums of dependent point processes to Poisson
processes.
\textit{Litovsk. Mat. Sb.} \textbf{15} 11--23, 223.
\MR{0418231}

\bibitem{Banys80}
Banys, R. (1980). On superpositions of random
measures and point processes.
In \textit{Mathematical Statistics and Probability Theory (Proc. Sixth Internat.
Conf., Wisla, 1978)} 26--37.
\textit{Lecture Notes in Statist.} \textbf{2}.
New York: Springer.
\MR{0577268}

\bibitem{Banis85}
Banys, R. (1985). A Poisson limit theorem for
rare events of a discrete random field.
\textit{Litovsk. Mat. Sb.} \textbf{25} 3--8.
\MR{0795851}



\bibitem{Barbour88}
Barbour, A.D. (1988). Stein's method and Poisson process
convergence. \textit{J. Appl. Probab.} \textbf{25(A)} 175--184.
\MR{0974580}

\bibitem{BB92}
Barbour, A.D. and Brown, T.C. (1992).
Stein's method and point process approximation. \textit{Stochastic
Process. Appl.} \textbf{43} 9--31.
\MR{1190904}

\bibitem{BE83} Barbour, A.D. and Eagleson, G.K. (1983).
Poisson approximation for some statistics based on
exchangeable trials. \textit{Adv. in Appl. Probab.} \textbf{15} 585--600.
\MR{0706618}

\bibitem{BHJ92}
Barbour, A.D., Holst, L. and Janson, S. (1992).
\textit{Poisson Approximation.} Oxford: Oxford Univ. Press.
\MR{1163825}

\bibitem{BarbourXia06}
Barbour, A.D. and Xia, A. (2006). Normal
approximation for random sums. \textit{Adv. in Appl. Probab.} \textbf{38}
693--728.
\MR{2256874}

\bibitem{Brown78}
Brown, T.C. (1978). A martingale approach
to the Poisson convergence of simple point processes.
\textit{Ann. Probab.} \textbf{6} 615--628.
\MR{0482991}

\bibitem{Brown79}
Brown, T.C. (1979).
Position dependent and stochastic thinning of point processes.
\textit{Stochastic Process. Appl.} \textbf{9} 189--193.
\MR{0548838}

\bibitem{BWX00}
Brown, T.C., Weinberg, G.V. and Xia, A. (2000).
Removing
logarithms from Poisson process error bounds.
\textit{Stochastic Process.
Appl.} \textbf{87} 149--165.
\MR{1751169}

\bibitem{BX95} Brown, T.C. and Xia, A. (1995).
On metrics in point process approximation.
\textit{Stochastics Stochastics Rep.} \textbf{52} 247--263.
\MR{1381671}


\bibitem{ChenShao04} Chen, L.H.Y. and Shao, Q.M. (2004). Normal
approximation under local dependence. \textit{Ann. Probab.} \textbf{32} 1985--2028.
\MR{2073183}

\bibitem{ChenXia04} Chen, L.H.Y. and Xia, A. (2004). Stein's
method, Palm theory and Poisson process
approximation. \textit{Ann. Probab.} \textbf{32} 2545--2569.
\MR{2078550}

\bibitem{Cinlar72}
\c{C}inlar, E. (1972). Superposition of point
processes. In
\textit{Stochastic Point Processes: Statistical Analysis, Theory, and
Applications (Conf., IBM Res. Center, Yorktown Heights, NY, 1971)} 549--606.
New York: Wiley.
\MR{0365697}


\bibitem{DaVerJon88} Daley, D.J. and Vere-Jones, D. (1988).
\textit{An Introduction to the Theory of Point Processes.}
New York: Springer.
\MR{0950166}


\bibitem{Fichtner75} Fichtner, K. (1975). Schwache Konvergenz von
unabh\"{a}ngigen \"{U}berlagerungen verd\"{u}nnte zuf\"{a}lliger
Punkfolgen. \textit{Math. Nachr.} \textbf{66} 333--341.
\MR{0370749}

\bibitem{Fichtner77} Fichtner, K. (1977).
Poissonsche zuf\"{a}llige Punktfolgen und ortsabh\"{a}ngige Verd\"
{u}nnungen. In
\textit{Transactions of the Seventh Prague Conference on Information Theory,
Statistical Decision Functions, Random Processes and of the Eighth
European Meeting of Statisticians (Tech. Univ. Prague, Prague, 1974)}
\textbf{A} 123--133. Dordrecht: Reidel.

\bibitem{Goldman67} Goldman, J.R. (1967). Stochastic point
processes: Limit theorems. \textit{Ann. Math. Statist.} \textbf{38} 771--779.
\MR{0217854}

\bibitem{GR96} Goldstein, L. and Rinott, Y. (1996). Multivariate
normal approximations by Stein's method and size bias couplings. \textit{J.
Appl. Probab.} \textbf{33} 1--17.
\MR{1371949}

\bibitem{GX06} Goldstein, L. and Xia, A. (2006). Zero biasing and
a discrete central limit theorem. \textit{Ann. Probab.} \textbf{34} 1782--1806.
\MR{2271482}

\bibitem{Grigelionis63} Grigelionis, B. (1963). On the convergence
of sums of random step processes to a Poisson process. \textit{Theory
Probab. Appl.} \textbf{8} 177--182.
\MR{0152013}


\bibitem{Jagers72} Jagers, P. (1972). On the weak convergence of
superpositions of point processes. \textit{Z. Wahrsch. Verw. Gebiete} \textbf{22} 1--7.
\MR{0309191}

\bibitem{Kallenberg75} Kallenberg, O. (1975).
Limits of compound and thinned point processes.
\textit{J. Appl. Probab.} \textbf{12} 269--278.
\MR{0391251}

\bibitem{Ka83} Kallenberg, O. (1983).
\textit{Random Measures.} London: Academic Press.
\MR{0818219}

\bibitem{Liese80} Liese, F. (1980).
\"{U}berlagerung verd\"{u}nnter und schwach abh\"{a}ngiger
Punktprozesse.
\textit{Math. Nachr.} \textbf{95} 177--186.
\MR{0592891}

\bibitem{Rachev91} Rachev, S.T. (1991). \textit{Probability Metrics and the
Stability of Stochastic Models.} Chichester: Wiley.
\MR{1105086}

\bibitem{Serfozo84} Serfozo, R. (1984). Thinning of cluster
processes: Convergence of sums of thinned point processes. \textit{Math. Oper.
Res.} \textbf{9} 522--533.
\MR{0769391}

\bibitem{Schuhmacher05a} Schuhmacher, D. (2005). Distance
estimates for dependent superpositions of point processes. \textit{Stochastic
Process. Appl.} \textbf{115} 1819--1837.
\MR{2172888}

\bibitem{Schuhmacher05b} Schuhmacher, D. (2005). Distance
estimates for Poisson process approximations of dependent thinnings.
\textit{Electron. J. Probab.} \textbf{10} 165--201 (electronic).
\MR{2120242}

\bibitem{Viswanathan92} Viswanathan, T. (1992). \textit{Telecommunication
Switching Systems and Networks.} New Delhi: Prentice-Hall of India.

\bibitem{Xia05} Xia, A. (2005). Stein's method and Poisson process
approximation. In \textit{An Introduction to Stein's Method}
(A.D. Barbour and L.H.Y. Chen, eds.) 115--181.
Singapore: World Scientific Press.
\MR{2235450}

\end{thebibliography}
\end{document}